\documentclass [12pt] {article}
\usepackage{amsfonts}
\newcommand{\qed} {\hspace {0.1in} \rule {1.5mm} {3.5mm} \vskip 0.1in}
\newcommand{\qedo} {\hspace {0.1in} \rule {1.5mm} {3.5mm}}

\newtheorem{lemma}{Lemma}[section]

\newtheorem{theorem}{Theorem}
\newtheorem{proposition}{Proposition}[section]
\newtheorem{definition}{Definition}[section]

\def\dc{\dim_\bC(W_n)}

\def\gi{g^{-1}}

\def\dim{{\rm dim}}

\def\<{\langle}
\def\>{\rangle}
\def\proof{\smallskip\noindent{\bf Proof:} }

\def\bF{{\mathbb F}}

\def\bC{{\mathbb C}}

\def\ca{\mbox{$\cal A$}}

\def\ch{\mbox{$\cal H$}}

\def\to{\rightarrow}
\def\ln{\lim_{n\to\infty}}

\title{The amenability and non-amenability of skew fields}
\author{{\sc G\'abor Elek}
\cr Mathematical Institute of
the Hungarian Academy of Sciences\cr P.O. Box 127, H-1364 Budapest, Hungary\cr
elek@renyi.hu}
\date{}
\begin{document}

\maketitle
\noindent{\bf Abstract.}  We investigate the amenability of skew field
extensions of the complex numbers. We prove that all skew fields of finite 
Gelfand-Kirillov transcendence degree are amenable. However there are both
amenable and non-amenable finitely generated skew fields of infinite Gelfand-Kirillov
transcendence degree.
\vskip 0.2in
\newpage
\section{Introduction}
\begin{definition}
Let $k$ be a commutative field and $\ca$ be a unital $k$-algebra. We say
that $\ca$ is amenable if for any finite subset
$\{r_1,r_2,\dots,r_n\}\subset \ca$ and real number $\epsilon>0$ there
exists a finite dimensional $k$- subspace $V\subset \ca$ such that
\begin{equation}
\frac{\dim_k (\sum^n_{i=1} r_iV)}{\dim_k V}<1+\epsilon
\end{equation}\end{definition}
In \cite{Elek}, we saw that the group algebra of an amenable group is
amenable, the group algebra of the free group of two generators is
non-amenable. Any affine algebra of subexponential growth or any
commutative algebra is amenable as well. It is easy to see that amenable
algebras satisfy the invariant basis number property. In this paper,
we investigate the amenability of skew field extensions. Our main result
is the following theorem.
\begin{theorem}

\noindent
\begin{description}
\item {(a)} If $k\subseteq D$ is a skew field that is finite dimensional over its center then
it is amenable.
\item {(b)} If $k\subseteq D$ has finite Gelfand-Kirillov transcendence degree
then it is amenable.
\item {(c)} There exist finitely generated amenable skew fields of
infinite Gelfand-Kirillov transcendence degree.
\item {(d)} If $k\subseteq D\subseteq E$ are skew fields and $D$ is non-amenable
then $E$ is non-amenable as well.
\item {(e)} The free field of Cohn \cite{Cohn} is non-amenable.
\end{description}
\end{theorem}
In the course of the proof we shall see that $\bC(\Gamma)$ is an amenable
$\bC$-algebra if and only if $\Gamma$ is an amenable group. The
corresponding  conjecture for general fields remains open.
\section{The Ore-property and amenability}
\begin{lemma}\label{3l}
If $\ca$ is an amenable domain, then it has the left Ore-property
that is for any non-zero $a,b\in\ca$; $a\ca\cap b\ca\neq 0$.
\end{lemma}
\proof
There exists a finite dimensional subspace $W\subseteq \ca$ such that
$\dim_k(aW\cap W)>\frac{1}{2}\dim_k (W)$ and
$\dim_k(bW\cap W)>\frac{1}{2}\dim_k (W)$, hence $aW\cap bW\neq 0$. \qed

\noindent
Consequently, if $\ca$ is an amenable domain, then one can consider
its classical ring of quotient $\widetilde{\ca}$, a skew field. 
\begin{proposition}\label{3p}
If $\ca$ is an amenable domain, then $\widetilde{\ca}$ is amenable as well.
\end{proposition}
\proof
Let $x_1,x_2,\dots, x_m\in \widetilde{\ca}$, then by the definition
of the classical ring of quotient there exists a non-zero element $r\in\ca$
such that $x_1r,x_2r,\dots, x_mr\in \ca$. By the amenability of $\ca$ there
exists a sequence of finite dimensional $k$-vector spaces $W_n\subseteq\ca$,
such that
\begin{equation}\label{3e}
\lim_{n\to\infty}
\frac{\dim_k(\sum^m_{i=1} x_ir W_n)}{\dim_k(W_n)}=1\,.
\end{equation}
\begin{lemma}\label{4l}
Let $Z_n=\{v\in W_n :\, x_irv\in W_n\,\mbox{for all $i$}\}$.
Then $\ln\frac{\dim_k(Z_n)}{\dim_k(W_n)}=1\,.$
\end{lemma}
\proof
By (\ref{3e}), for any $1\leq i \leq m$, $\ln\frac{\dim_k(x_irW_n\cap W_n)}
{\dim_k(W_n)}=1.$ Hence
if \\
 $Z_n^i=\{v\in W_n,\, x_irv\in W_n\}$, then $\ln
\frac{\dim_k(Z_n^i)}{\dim_k(W_n)}=1\,.$
Since $Z_n=\cap^m_{i=1} Z_n^i$, the lemma follows. \qed

\noindent
Obviously, if $T_n=\{v\in W_n:\, rv\in W_n\}$ then
$\ln \frac{\dim_k(T_n)}{\dim_k(W_n)}=1\,.$ Thus if $S_n=T_n\cap Z_n$, then
$\ln \frac{\dim_k(S_n)}{\dim_k(W_n)}=1\,$, and 
$\ln \frac{\dim_k(rS_n)}{\dim_k(W_n)}=1\,$.
Now if $v\in rS_n$, then $v\in W_n$ and for any $1\leq i \leq m$, $x_iv\in
W_n$,
that is $rS_n\subseteq W_n$ and $ x_irS_n\subseteq W_n$.

\noindent
Since $\ln\frac{\dim_k(W_n)-\dim_k(rS_n)}{\dim_k(W_n)}=0$, it follows that
$$\ln\frac{\dim_k(\sum^m_{i=1} x_i W_n)}{\dim_k(W_n)}=1\,,$$
proving the amenability of $\widetilde{\ca}$. \qed

\noindent
Notice that by our previous proposition the classical ring of fraction
of the Weyl-algebras are amenable. Let us recall the notion of
the Gelfand-Kirillov transcendence degree \cite{GK}.
\begin{definition}
Let $D$ be a skew field extension of the commutative field $k$.
Then
$$ \mbox{GK-tr deg}(D)=\sup_V \inf_{0\neq r} \limsup_{n\to\infty}
\frac{\log \dim_k((k+Vr)^n)}{\log n}\,,$$
where $V$ runs through all the finite
dimensional subspaces of $D$ containing the unit.
\end{definition} Note that
$\mbox{GK-tr deg}(D)<\infty$ if and only if there exists a $d>0$ with the
following property: 
For any $c_1,c_2,\dots c_m\in D$
there exists $r\in D$ such that
the affine algebra generated by $1, c_1r,c_2r,\dots c_mr$ has polynomial growth
with Gelfand-Kirillov dimension not greater then $d$.
\begin{proposition}\label{5p}
If $D$ has finite Gelfand-Kirillov transcendence degree then 
\\  $D$ is amenable.
\end{proposition}
\proof
We prove that if for any subset $c_1,c_2,\dots,c_m\in D$ there exists
$r\in D$ such that $1, c_1r, c_2r,\dots, c_mr$ generate an amenable domain,
then $D$ is amenable. Indeed, if
$$\ln \frac {\dim_k(\sum_{i=1}^mc_irW_n)}{\dim_k(W_n)}=1\,,$$
then 
$$\ln \frac {\dim_k(\sum_{i=1}^mc_iV_n)}{\dim_k(V_n)}=1\,,$$
where $V_n=rW_n$. Thus $D$ is in fact amenable. \qed

\noindent
Note however,
that Aizenbud proved in \cite{Aiz} that for some poly-cyclic groups $\Gamma$
the classical ring of quotient of the
group algebra $k(\Gamma)$ has infinite Gelfand-Kirillov transcendence
degree. Since such groups are amenable, these finitely generated skew fields above are amenable
as well.

\noindent 
Also, since all affine commutative algebras and all matrix algebras
over such algebras are of polynomial growth, we have the following
corollary:
\begin{proposition} \label{6p1}
If $k\subseteq L$ is finite dimensional over its center,
then $L$ is amenable.
\end{proposition}
The class of amenable skew fields are closed under some operations:
\begin{proposition} \label{elso}

\noindent
\begin{description}
\item{(a)} The direct limit of amenable skew fields is amenable.
\item{(b)} If $k\subseteq D$ and any finitely generated subskewfield
$E\subseteq D$ is amenable, then $D$ is amenable
as well.
\item{(c)} If $A$ and $B$ are amenable skew fields and $A\otimes_k B$ is
a domain, then $A\otimes_k B$ is amenable, hence by Proposition \ref{3p}
$\widetilde{A\otimes_k B}$ is also an amenable skew field.
\end{description}
\end{proposition}
\proof
The parts (a) and (b) are easy exercises, for (c), let us consider the
amenable $k$-algebras $A$ and $B$. It is enough to prove that if
$a_1,a_2,\dots,a_l\in\ca$, $b_1,b_2,\dots,b_m\in B$ then there exists
a sequence of finite dimensional subspaces $\{V_n\}\subseteq A\otimes_k B$
such that
$$\ln \frac{\dim_k(\sum^l_{i=1}\sum^m_{j=1} (a_i\otimes b_j)V_n)}{\dim_k (V_n)}=1\,.$$
Let $\{W_n\}\subseteq A, \{Z_n\}\subseteq B$ be sequences
of finite dimensional subspaces such that
$$\ln \frac{\dim_k(\sum^l_{i=1}a_iW_n)}{\dim_k (W_n)}=1\,,$$
and
$$\ln \frac{\dim_k(\sum^m_{j=1}b_jZ_n)}{\dim_k (Z_n)}=1\,.$$
Set $V_n=W_n\otimes Z_n$. Then
$$(\sum^l_{i=1}\sum^m_{j=1} (a_i\otimes b_j)V_n)\subseteq
(\sum^l_{i=1}a_iW_n)\otimes
(\sum^m_{j=1}b_jZ_n)\,,$$
hence the proposition follows. \qed
\section{Non-amenable subskewfields}
\begin{proposition}
Let $k\subseteq D\subseteq E$ be skewfields. Suppose that $D$ is non-amenable,
then $E$ is non-amenable as well.
\end{proposition}
\proof
If $D$ is non-amenable, then there exists $\alpha>1$ and 
$g_1,g_2,\dots,g_l\in D$ such that for any finite dimensional
$k$-subspace $W\subseteq D$ 
\begin{equation} \label{e81}
\dim_k(\sum^l_{i=1}g_iW)\geq \alpha\, \dim_k (W)\quad.
\end{equation}
It is enough to prove that for any $m\geq 1$ and any finite dimensional
$k$-subspace $Z$ of the $m$-dimensional left $D$-module $\oplus^m_{j=1} D$
\begin{equation} \label{e82}
\dim_k(\sum^l_{i=1}g_iZ)\geq \alpha\, \dim_k (Z)\quad.
\end{equation}
We proceed by induction. The $m=1$ case is just (\ref{e81}).
Suppose that (\ref{e82}) holds for $m=1,2,\dots n$. Let $Z\subseteq
\oplus^{n+1}_{j=1} D$ and $\pi^{n+1}: \oplus^{n+1}_{j=1} D\to D$ the
projection onto the last coordinate. \\
If $\dim_k (\pi^{n+1} (Z))=\dim_k (Z)$, then
$$\dim_k(\sum^l_{j=1}g_jZ)=\dim_k(\sum^l_{j=1} g_j(\pi^{n+1}(Z)))\geq
\alpha\,\dim_k(Z)\,.$$
If $\dim_k(\pi^{n+1}(Z))<\dim_k(Z)$, then
$Z=Z_1\oplus Z_2$ as $k$-subspaces, where
$Z_1\subseteq \oplus^n_{j=1} D$ and $\dim_k(\pi^{n+1}(Z))=\dim_k(Z_2)$.
Let $v_1, v_2,\dots, v_{s_1}$ be a $k$-basis for $Z_1$ and
let $w_1, w_2, \dots, w_{s_2}$ be a $k$-basis for $Z_2$. Let $G$ be the
$k$-vectorspace generated by $\{g_1,g_2,\dots,g_l\}$.
 Then by the
inductional hypothesis, there
exists elements 
$$g^i_j\in G, 1\leq i \leq t_1, 1\leq j \leq s_1$$
$$h^i_j\in G, 1\leq i \leq t_2, 1\leq j \leq s_2\,,$$
such that $t_1\geq \alpha\,s_1$, $t_2\geq \alpha\,s_2$ and :
$$\{z_i=\sum_{j=1}^{s_1} g^i_j v_j\}_{i=1,2,\dots,t_1}$$
and
$$\{z'_i=\sum_{j=1}^{s_2} h^i_j \pi^{n+1}(w_j)\}_{i=1,2,\dots,t_2}$$
are both independent system of vectors in $\oplus^n_{i=1} D$ resp.
$\pi^{n+1}(\oplus^{n+1}_{i=1} D)$.

\noindent
Then $\{z_i=\sum_{j=1}^{s_1} g^i_j v_j\}\cup
 \{z'_i=\sum_{j=1}^{s_2} h^i_j (w_j)\}$ form an independent
system of more than $\alpha\cdot\dim_k(Z)$ elements, proving (\ref{e82})
for the subspace $Z\subseteq \oplus^{n+1}_{j=1} D$ \qed
\section{The construction of a non-amenable skew field}
\begin{definition}
Let $\Gamma$ be a discrete group and $V$ be a vectorspace over $k$.
A representation $\phi:\Gamma\to GL(V)$ is called algebrically
amenable, if for any
$g_1, g_2,\dots g_m\in\Gamma$ and $\epsilon>0$ there exists a finite
dimensional vectorspace $W\subseteq V$ such that
$$\frac{\dim_k(\sum^m_{i=1} g_iW)}{\dim_k (W)}\leq 1+\epsilon\,.$$
\end{definition}
Note that if $\ca$ is an amenable $k$-algebra, then the group
of invertible elements of $\ca$ acts on $\ca$ in an algebrically
 amenable fashion. We should observe that the group of non-zero
elements of an amenable skew field can be in fact a  non-amenable group.
As a matter of fact it is Lichtman's conjecture that any skew field
which is infinite dimensional over its center contains a free group
of two generators. By the result of Makar-Limanov \cite{Mak} the
first Weyl skew field contains even a free group algebra
as a subalgebra. 
\begin{proposition}\label{bekka}
Let $\Gamma$ be a countable group represented by
unitary transformations on a Hilbert-space $\ch$.
If the representation is algebraically amenable, then
it is amenable in the sense of Bekka \cite{Bek} that is
there exists a linear map $\tau:B(\ch)\to\bC$ from the
space of bounded linear operators to the complex numbers
such that $\tau(Id)=1$ and
$$\tau(A)=\tau(\gi A g)\,,$$
for any $A\in B(\ch)$, $g\in G$.
\end{proposition}
\proof
By our assumption, there exists a sequence of
finite dimensional subspaces $W_n\subseteq\ch,\,n\geq 1$ such that for
any $g\in \Gamma$:
$$\ln\frac{\dim_\bC(gW_n+W_n)}{\dim_\bC(W_n)}=1\,.$$
Let $P_n$ denote the orthogonal projection onto $W_n$. Then for any
$A\in B(\ch)$,
$$\frac{Tr(P_nAP_n)}{\dc}$$
is a sequence of complex number bounded by the norm of $A$.
Hence one can consider
$$\tau(A)=\lim_\omega\frac{Tr(P_nAP_n)}{\dc}\quad,$$
where $\omega$ is an ultrafilter on the natural numbers and
$\lim_{\omega}$ is the corresponding ultralimit. Then the map
$\tau:B(\ch)\to\bC$ is clearly bounded and linear.
In order to prove that
$\tau(A)=\tau(gA\gi)$ it is enough to see that for any $g\in\Gamma$ and $A\in
B(\ch)$:
$$\ln\frac{Tr(P_nAP_n)-Tr(P_n\gi Ag P_n)}{\dc}=0.$$
Let $G_n=\gi(gW_n\cap W_n)$. Then
$\ln\frac{\dim_\bC(G_n)}{\dc}=1$ and both $G_n$ and $gG_n$ are subspaces
of $W_n$. Pick an orthonormal basis $\{v_1,v_2,\dots,v_{k_n}\}$ for $G_n$ and extent it to
an orthonormal basis for the whole $W_n$ by adding
 $w_1,w_2,\dots,w_{l_n}$, where
$k_n=\dim_\bC(G_n),\,l_n=\dim_\bC(W_n)-\dim_\bC(G_n)$.

\noindent
Let $g_n$ be a unitary transformation on $W_n$ such that $g_n\mid_{G_n}=g$.
Define $g_n$ to be the identity transformation on $W_n^\perp$. Then:
\begin{equation}
\label{e12}
Tr(P_n\gi_nAg_n P_n)=Tr(P_nAP_n)\quad.
\end{equation}
Indeed, $g_n$ commutes with $P_n$, hence 
$$Tr(P_n\gi_nAg_n P_n)=Tr(\gi_n P_nAP_ng_n)$$
and
$$Tr(\gi_nP_nAP_ng_n)=Tr(P_nAP_n)$$
by the property of the trace.
\begin{lemma}\label{12l}
$$\ln\frac{Tr(P_n\gi_n Ag_n P_n)-Tr(P_n\gi AgP_n)}{\dc}=0\,.$$
\end{lemma}
\proof Note that:
$$Tr(P_n\gi AgP_n)=\sum^{k_n}_{i=1}\langle\gi A g(v_i),v_i\rangle+
\sum^{l_n}_{i=1}\langle\gi A g(w_i),w_i\rangle= $$
$$=\sum^{k_n}_{i=1}\langle Ag(v_i),g(v_i)\rangle+\sum^{l_n}_{i=1}\langle
Ag(w_i),g(w_i)\rangle\,.$$
Similarly,
$$Tr(P_n\gi_n  Ag_n P_n)=
\sum^{k_n}_{i=1}\langle Ag_n(v_i),g_n(v_i)\rangle+\sum^{l_n}_{i=1}\langle
Ag_n(w_i),g_n(w_i)\rangle\quad.$$
Since $g_n(v_i)=g(v_i)$ we have the following estimate:
$$|Tr(P_n\gi_n Ag_n P_n)-Tr(P_n\gi AgP_n)|\leq $$
$$\sum^{l_n}_{i=1}(|\langle
Ag_n(w_i),g_n(w_i)\rangle|+|\langle
Ag(w_i),g(w_i)\rangle|\leq 2\|A\|\cdot l_n\quad.$$
Therefore the lemma follows from the fact that
$$\ln\frac{l_n}{\dc}=0\,.\quad\quad\qedo$$
Consequently, by the Lemma \ref{12l} and (\ref{e12}) we can see that
$\tau(A)=\tau(\gi A g)$.\qed

\begin{proposition}\label{p13}
If $\Gamma$ is a countable group, then $\bC(\Gamma)$ is amenable as a
$\bC$-algebra if and only if $\Gamma$ is an amenable group. Moreover
$\Gamma$ acts on $l^2(\Gamma)$ in an algebrically amenable fashion.
\end{proposition}
\proof
If $\Gamma$  amenable then $k(\Gamma)$ is amenable for any field $k$
\cite{Elek}. Now, let us suppose that $\Gamma$ is non-amenable.
 If $\bC(\Gamma)$
were amenable, then the natural representation of $\Gamma$ on $\bC(\Gamma)$
would be algebrically amenable. Therefore the natural representation of $\Gamma$
by translations on the Hilbert-space $\bC(\Gamma)$ must be algebrically
amenable as well. By our previous proposition, there exists a bounded
linear functional $\tau:B(l^2(\Gamma))\to\bC$ such that
$\tau(A)=\tau(\gi A g)\,,$
for any $A\in B(\ch)$, $g\in \Gamma$. Consider the commutative $C^*$-algebra
$l^\infty(\Gamma)\subseteq B(l^2(\Gamma))$ acting on $l^2(\Gamma)$
by pointwise multiplication. That is if $f\in l^\infty(\Gamma),\,G\in
l^2(\Gamma)$:
$$T_f(G)(\delta)=f(\delta)G(\delta)\,.$$
The group $\Gamma$ acts on both $l^\infty(\Gamma)$ and $l^2(\Gamma)$ by
translations:
$$(\gamma f)(\delta)=f(\gamma^{-1}\delta)\,.$$
\begin{lemma}
$$T_{\gamma f}=\gamma\circ T_f\circ \gamma^{-1}\,.$$
\end{lemma}
\proof
Let $1_\delta\in l^2(\Gamma)$ be the function which vanishes everywhere
except at $\delta$ where it takes the value $1$. 

\noindent
Then $\gamma^{-1}(1_\delta)=1_{\gamma^{-1}\delta}$ and
$(T_f\circ\gamma^{-1})(1_\delta)=f(\gamma^{-1}\delta)\cdot
1_{\gamma^{-1}\delta}$.
Thus $(\gamma\circ T_f\circ
\gamma^{-1})(1_\delta)=f(\gamma^{-1}\delta)1_\delta\,.$

\noindent
On the other hand, $T_{\gamma\,f}(1_\delta)=(\gamma f)(\delta)\cdot 1_\delta=
f(\gamma^{-1}\delta)1_\delta\,.$\qed

\noindent
Hence, $\tau$ would define a translation invariant
linear functional on $l^\infty(\Gamma)$, which takes the value $1$
on the constant one function. This is in contradiction with the
non-amenability of the group $\Gamma$. Thus our proposition holds.\quad\qed
\vskip 0.2in
\noindent
{\bf Remark:} At the 2003 Gaeta Conference, Zelmanov mentioned the following
conjecture:
If a discrete group of Kazhdan's property $(T)$ is represented on an infinite 
dimensional vectorspace without a finite dimensional invariant subspace, then
the representation is not algebraically amenable.

\noindent
If the action in the conjecture is in fact unitary, then by the result
of \cite{BV}, the representation can not be amenable in the sense
of Bekka. Hence by Proposition \ref{bekka}, Zelmanov's conjecture holds for
unitary representations.

\vskip 0.15in
\noindent
Now, we are in the position to construct the non-amenable skew field.
Let us consider the von Neumann algebra $W(\bF_2)$ of the free group
of two generators. The algebra $W(\bF_2)$ satisfies the Ore-condition
with respect to its non-zero divisors and thus it imbeds to the ring of 
affiliated operators $U(\bF_2)$ By Lemma 10.51 of \cite{Lueck}
there exists a skew field $D(\bF_2)$ such that
$\bC(\bF_2)\subseteq D(\bF_2)\subseteq U(\bF_2)$ namely the division closure of $\bC(\bF_2)$
in $U(\bF_2)$.

\begin{proposition}\label{utolso}
$D(\bF_2)$ is a non-amenable skew field.
\end{proposition}
\proof
Recall \cite{Lueck} that the von Neumann algebra can be identified
with the dense linear subspace of $l^2(\bF_2)$ consisting of
those vectors $w\in l^2(\bF_2)$ such that
the convolution by $w$ defines a bounded operator on $l^2(\Gamma)$. 
Thus the natural
action of the group $\bF_2$ on the vectorspace $W(\bF_2)$ must be
non-amenable. 
Since $\bC(\bF_2)\subseteq D(\bF_2)\subseteq U(\bF_2)$, it is
enough to prove that the action of $\bF_2$ on $U(\bF_2)$
is not algebrically amenable.
Let $V_n\subseteq U(\bF_2)$ be a sequence of linear subspaces
such that for any $g\in \bF_2$
$$\ln\frac{\dim_\bC(gV_n+V_n)}{\dim_\bC(V_n)}=1\,.$$
Then by the Ore-property, there exists non-zero divisors
$s_n\in W(\bF_2)$ such that
$V_ns_n\subseteq W(\bF_2)$. Then
$$\ln\frac{\dim_\bC(gV_ns_n+V_ns_n)}{\dim_\bC(V_ns_n)}=1\,.$$
This would mean the  action of $\bF_2$ on $U(\bF_2)$ is algebrically
amenable, leading to a contradiction. \qed

\noindent
Linnell \cite{Lin} proved that $D(\bF_2)$ is the free field.
This finishes the proof of our Theorem 1. Note that Proposition \ref{utolso}
and Proposition \ref{elso} imply the existence of finitely generated
non-amenable skew fields.

\end{document}